\documentstyle[12pt]{article}
\newcommand{\n}{\noindent}

 \newcommand{\be}{\begin{equation}}
 \newcommand{\ee}{\end{equation}}
 
 \newfont{\gl}{eufm10 scaled \magstep1}  % german letters as well
 
 \newcommand{\bC}{{\bf C}}
 
\newcommand{\bH}{{\bf H}}
 \newcommand{\bI}{{\bf I}}
 
 \newcommand{\bP}{{\bf P}}
 \newcommand{\bR}{{\bf R}}

 \newcommand{\bZ}{{\bf Z}}
 \newcommand{\cA}{{\cal A}}
 \newcommand{\cB}{{\cal B}}
 \newcommand{\cC}{{\cal C}}
 
 \newcommand{\cE}{{\cal E}}
 \newcommand{\cF}{{\cal F}}
 \newcommand{\cG}{{\cal G}}
 \newcommand{\cH}{{\cal H}}

 \newcommand{\cL}{{\cal L}}

 \newcommand{\cO}{{\cal O}}

 \newcommand{\cV}{{\cal V}}

 \newcommand{\glu}{\hbox{\gl u}} % gothic u
 \newcommand{\gls}{\hbox{\gl s}} % gothic s
 \newcommand{\ra}{\rightarrow}
 \newcommand{\lra}{\longrightarrow}
 
 \newcommand{\kahler}{K\"{a}hler}

 \newcommand{\mod}{\mathop{{\fam0 mod}}\nolimits}
 \newcommand{\End}{\mathop{{\fam0 End}}\nolimits}
  \newcommand{\ad}{\mathop{{\fam0 ad}}\nolimits}
 \newcommand{\rank}{\mathop{{\fam0 rank}}\nolimits}

 \newcommand{\tr}{\mathop{{\fam0 Tr}}\nolimits}
 \newcommand{\Hom}{\mathop{{\fam0 Hom}}\nolimits}

 \newcommand{\vol}{\mathop{{\fam0 Vol}}\nolimits}
 \newcommand{\pic}{\mathop{{\fam0 Pic}}\nolimits}
 \newcommand{\U}{\mathop{{\fam0 U}}\nolimits}
 \newcommand{\so}{\mathop{{\fam0 SO}}\nolimits}
 \newcommand{\spin}{\mathop{{\fam0 Spin}}\nolimits}
 \newcommand{\spinc}{\spin^c}
\newcommand{\ch}{\mathop{{\fam0 ch}}\nolimits}
 
 \newcommand{\YMH}{\mathop{{\fam0 YMH}}\nolimits}
\newcommand{\SW}{\mathop{{\fam0 SW}}\nolimits}
 
 \newcommand{\lie}{\mathop{{\fam0 Lie }}\nolimits}
 
 \newcommand{\dbar}{\overline{\partial}}
 \newtheorem{definition}{Definition}[section]
 \newtheorem{lemma}[definition]{Lemma}
 \newtheorem{prop}[definition]{Proposition}
 \newtheorem{thm}[definition]{Theorem}

 \newcommand{\pf}{{\em Proof}. }
 \newcommand{\remark}{\noindent{\em Remark}. }
 
 \newcommand{\ps}{{p^\ast}}
 \newcommand{\qs}{{q^\ast}}
 \newcommand{\su}{{SU}}
\newcommand{\xp}{{{X\times\bP^1}}}

 \newcommand{\hA}{{\hat{A}}}
 \newcommand{\hb}{{\hat{b}}}
 
 \newcommand{\hL}{{\hat{L}}}
 \newcommand{\hE}{{\hat{E}}}

 \newcommand{\cves}{coupled vortex equations}
 
 \newcommand{\ts}{$\tau$-stable}

 \newcommand{\ymh}{Yang--Mills--Higgs}
 \newcommand{\he}{Hermitian--Einstein}

 \newcommand{\JDG}{J. Diff. Geom.}
 \newcommand{\qed}{\hfill$\Box$}
 
 \newcommand{\hk}{Hitchin--Kobayashi}
 \newcommand{\sw}{Seiberg--Witten}

 \newcommand{\extn}{{0\lra \cE_1\lra \cE\lra  \cE_2\lra 0}}
 \newcommand{\sextn}{{0\lra \cE_1'\lra\cE'\lra  \cE_2'\lra 0}}

\newcommand{\slp}{{S_L^+}}
\newcommand{\slm}{{S_L^-}}
\newcommand{\slpm}{{S_L^\pm}}

\begin{document}

\title{Non-abelian monopoles and vortices}
\author{Steven B. Bradlow \and Oscar Garc\'{\i}a--Prada}

\date{}

\maketitle

%%%%%%%%%%%%%%%%%%%%%%%%%%%%%%%%%%%%%%%%%%
\section{Introduction}\label{introduction}
%%%%%%%%%%%%%%%%%%%%%%%%%%%%%%%%%%%%%%%%%%

In the short time since their discovery, the Seiberg--Witten equations
have already proved to be a powerful tool in the study of smooth
four-manifolds.
Virtually all the hard-won gains that have been obtained using the heavy
machinery of
Donaldson invariants, can be recovered with a fraction of the effort
if the (SU(2)) anti-self-duality equations are replaced by the
($\U(1)$) Seiberg--Witten equations. In addition, the new equations
probe deep features of symplectic structures.  They have also been used
to study geometric questions on K\"ahler surfaces. (See \cite{D3} for a
useful  survey.)

The  impressive success of the original equations has naturally led
to speculation about possible generalizations and other related sets
of equations.  The  original equations as proposed by Seiberg--Witten
are associated with a Hermitian line bundle, and thus with the
abelian group $\U(1)$.  One way to generalize the equations is thus
to look for versions based on larger, non-abelian groups.  This
means replacing the line bundle with a higher rank complex vector
bundle.  Indeed a number of authors have proposed versions of the
equations along these lines. These include, among others,
Okonek and Teleman \cite{OT1,OT2},
Pidstrigach and Tyurin \cite{PT},  Labastida and Mari\~no \cite{LM},
as well as the  second author \cite{G5}.  Some of these (cf. \cite{PT})
play a key role in attempts to prove the conjecture of Witten \cite{W}
concerning the equivalence of the old Donaldson and the new
\sw\ invariants.
(See also \cite{D3}, and \cite{FL} for more recent progress in
this direction.)

It is striking that no two of the above mentioned authors
consider precisely the same set of equations. One conclusion to be drawn from
this abundance of equations, is
that there is apparently more than one natural way to write down non
abelian versions of the Seiberg-Witten equations.  This leads to the
question: Are  some versions more reasonable, or more natural, than others?
The material in this paper gives one perspective on this question.

The main idea in our point of view is to exploit the special form of
the Seiberg-Witten equations in the case where the four manifold is
a K\"ahler surface.  In this case the original Seiberg--Witten
equations are known to reduce essentially to familiar equations in gauge
theory known as the abelian vortex equations.  Looked at from
the opposite direction, the Seiberg--Witten equations serve as
a ''Riemannian version'' of the vortex equations.  The key point is
that there are a number of well motivated, natural generalizations
of such vortex equations.  All of these are defined on complex
vector bundles over K\"ahler manifolds, and thus in particular over
K\"ahler surfaces.  Our guiding principle is
that the generalizations of the Seiberg--Witten equations should provide
''Riemannian versions'' of these vortex-type equations over K\"ahler
surfaces.

In this paper we explore essentially two such non-abelian generalizations.
We also make some remarks concerning a different aspect of the
relation between the vortex and the Seiberg-Witten equations.
This aspect has to do with the parameters which appear in the
vortex equations.  In their
original form, these were taken to be real numbers, i.e. constant
functions on the base manifolds. In the versions that emerge from the
Seiberg-Witten equations, the analogous terms turn out to be non-constant
functions (related to the scalar curvature).  This has prompted a closer
look at the affected terms in the vortex equation.   We discuss various ways
of
incorporating ---  and interpreting --- this level of generality in the
analysis  of the vortex equations.

In the interests of completeness, we have included a certain amount of
standard background material on the Seiberg-Witten and vortex equations.
%%%%%%%%%%%%%%%%%%%%%%%%%%%%%%%%%%%%%%%%%%

\noindent {\bf Acknowledgements.}
The second author wishes to thank the organisers of the Aarhus
Conference in Geometry and Physics, and especially J{\o}rgen Andersen,
for their kind invitation to participate in the Conference,
and to visit the Mathematics Institute in  Aarhus,
as well as for their warm hospitality.

%%%%%%%%%%%%%%%%%%%%%%%%%%%%%%%%%%%%%%%%%%%%%%%%%
\section{The \sw\ monopole equations}\label{sw}
%%%%%%%%%%%%%%%%%%%%%%%%%%%%%%%%%%%%%%%%%%%%%%%%%
%%%%%!!!!!!!! CHANGES STARTING HERE
In this section we briefly review the \sw\
equations and the analysis of these equations in the \kahler\ case.
For more details, see the original papers by Witten \cite{W}
and Kronheimer and Mrowka \cite{KM}, or any recent survey on the subject
(e.g. \cite{D3,G5}).

Let  $(X,g)$ be a compact, oriented, Riemannian four-manifold.
To write the \sw\ equations one needs a $\spinc$-structure on $X$.
This involves the choice of a Hermitian line bundle $L$ on $X$
satisfying that $c_1 (L)\equiv w_2(X)\;\mod\;2$.
A $\spinc$-structure is then a lift of the fibre product of the
$\so(4)$-bundle of orthonormal
frames of $(X,g)$ with the $\U(1)$-bundle defined by $L$ to a
$\spinc(4)$-bundle, according to the short exact sequence
$$
0\lra\bZ_2\lra\spinc(4)\lra\so(4)\times \U(1)\lra 1.
$$
Using the two fundamental irreducible 2-dimensional representations
of $\spinc(4)$---the so-called
Spin representations---we can construct the associated  vector
bundles of positive and negative spinors $S_L^\pm$.
These are rank 2 Hermitian vector bundles whose determinant is $L$
\cite{H,LaMi}.
 The set of $\spinc$-structures on $X$
is thus parametrised, up to the finite group $H^1(X,\bZ_2)$, by
$$
\spinc(X)=\{c\in H^2(X,\bZ)\;|\; c\equiv w_2(X)\; \mod\;2\}.
$$
%%%%%%% ENDING HERE

Let us fix a $\spinc$-structure $c\in \spinc(X)$, and let $L=L_c$
be the corresponding Hermitian line bundle,  and $S_L^\pm$ the
corresponding spinor bundles.
The \sw\ {\em monopole equations} are equations for a pair $(A,\Psi)$
consisting of a unitary connection on $L$ and a smooth section of $\slp$.
Using the connection $A$ one has  the Dirac operator
$$
D_A:\,\Gamma(\slp)\lra\Gamma(\slm).
$$
The first condition is that $\Psi$ must be in the kernel of the
Dirac operator.

The curvature $F_A\in\Omega^2=\Omega^+\oplus\Omega^-$,
can be decomposed in the self-dual and anti-selfdual parts
$$
F_A=F_A^+ + F_A^-.
$$
Using  the spinor $\Psi$ we can consider another self-dual
2-form that we may couple to $F_A^+$ to obtain our second equation.
Let $\ad_0\slp$ be the subbundle of the adjoint bundle of $\slp$
consistig of the traceless skew-Hermitian endomorphisms --- its fibres
are hence isomorphic to $\gls\glu(2)$.
We have a map
$$
\Omega^0(\slp)\lra\Omega^0(\ad_0\slp)
$$
$$
\Psi\mapsto i(\Psi\otimes\Psi^\ast)_0\ ,
$$
where $\Psi^\ast$ is the adjoint of $\Psi$, and the 0 subindex means that
we are taking the trace-free part.
This map is fibrewise modelled on the map $\bC^2\lra\gls\glu(2)$,
given by $v\mapsto i(v\overline{v}^t)_0$.

One of the basic ingredients  that makes the \sw\ equations possible
is the identification between the space of self-dual 2-forms and the
skew-Hermitian automorphism of the positive spin representation
\cite{AHS}.
 This is a basic fact in Clifford algebras in dimension four, that takes
place at each point of the manifold, and that can be carried out over
the whole manifold precisely when one has a $\spinc$-structure.
More specifically, we have the isomorphism
\be
\ad_0\slp\cong\Lambda^+\ . \label{iso}
\ee
We can now interpret
$i(\Psi\otimes\Psi^\ast)_0$ as a section of $\Lambda^+$, i.e. as an
element in $\Omega^+$.
The monopole equations consist in the system of equations
\be
\left.\begin{array}{l}
D_A \Psi=0\\
F_A^+=i(\Psi\otimes\Psi^\ast)_0
\end{array}\right \}. \label{me}
\ee
In writing the second equation there is an  abuse of notation, since we
are not especifying what the isomorphism (\ref{iso}) is. Notice also that
$F_A^+$ is a purely imaginary self-dual 2-form, and hence we are in fact
identifying $\ad_0\slp$ with $i\Lambda^+$.
%%%%%%%%%%%%%%%%%%%%%%%%%%%%%%
%\subsection{The \kahler\ case}
%%%%%%%%%%%%%%%%%%%%%%%%%%%%%%

We shall analyse now the monopole equations in the case in
which $(X,g)$ is \kahler.
Recall that a \kahler\ manifold is Spin if and only if there exists
a square root of the canonical bundle $K^{1/2}$ \cite{A,H}.
Moreover the spinor bundles are
$$
S^+=(\Lambda^0\oplus\Lambda^{0,2})\otimes K^{1/2}=K^{1/2}\oplus K^{-1/2}
$$
$$
S^-=\Lambda^{0,1}\otimes K^{1/2}.
$$
In this situation the spinor bundles for the $\spinc$-structure $c$ are
given by $S_L^\pm=S^\pm\otimes L^{1/2}$ (notice that
$L^{1/2}$ exists since $c_1(L)\equiv 0\;\mod\;2$).
Even if $X$ is not Spin, i.e. even if $K^{1/2}$ and
$L^{1/2}$ do not exist, the bundles $S_L^\pm$ do exist. In other words,
there exists a square root of $K\otimes L$. Let us denote
$$
\hL=(K\otimes L)^{1/2}.
$$
Then
$$
\slp=\hL\oplus \Lambda^{0,2}\otimes \hL,\;\;\slm=\Lambda^{0,1}\otimes\hL
$$
and
$$
\Gamma(\slp)=\Omega^0(\hL)\oplus\Omega^{0,2}(\hL).
$$
We can write $\Psi$ according to this decomposition as a pair
$\Psi=(\phi,\beta)$. The Dirac operator can be written in this language
as
$$
\dbar_\hA +\dbar_\hA^\ast\;\;:
\Omega^0(\hL)\oplus\Omega^{0,2}(\hL)\lra\Omega^{0,1}(\hL),
$$
where $\dbar_\hA$ is the $\dbar$ operator on $\hL$ corresponding to the
connection $\hA$ on $\hL$ defined by the connection $A$ on $L$ and the
metric connection on $K$ (cf. \cite{H}).

On the other hand recall that
$$
\Lambda^+\otimes\bC=\Lambda^0\omega\oplus\Lambda^{2,0}
\oplus\Lambda^{0,2},
$$
where $\omega$ is the \kahler\ form.
According to this decomposition, the isomorphism (\ref{iso})
(or rather $\ad_0\slp\cong i\Lambda^+$)
is explicitely given by
\be
i(\Psi\otimes\Psi^\ast)_0\mapsto i(|\phi|^2-|\beta|^2)\omega +
\beta\overline{\phi}-\phi\overline{\beta}.
\ee
We may thus write the monopole equations (\ref{me}) as
\be
 \begin{array}{l}
\dbar_\hA\phi+\dbar^\ast_\hA\beta=0\\
\Lambda F_A=i(|\phi|^2-|\beta|^2) \\
F_A^{2,0}=-\phi \overline{\beta}\\
F_A^{0,2}=\beta\overline{\phi}
\end{array}\label{kme}
\ee
where $\Lambda F_A$ is the contraction of the curvature with the
\kahler\ form.

It is not difficult to see (cf. \cite{W}) that the solutions to these
equations are such that either $\beta=0$ or $\phi=0$, and moreover it
is not
possible to have irreducible solutions, i.e. solutions with $\Psi\neq 0$,
of both types simultaneously for a fixed $\spinc$-structure.
We thus have one of the following two situations:

(i)\ $\beta=0$ and the equations reduce to
\be
\begin{array}{l}
F_A^{0,2}=0\\
\dbar_\hA\phi=0\\
\Lambda F_A=i|\phi|^2
\end{array}\label{ve1}
\ee

(ii)\ $\phi=0$ and then
\be
\begin{array}{l}
F_A^{0,2}=0\\
\dbar_\hA^\ast\beta=0\\
\Lambda F_A=-i|\beta|^2\ .
\end{array} \label{ve2}
\ee
\remark
We have omitted the equation $F_A^{2,0}=0$, since by unitarity of the
connection this is equivalent to $F_A^{0,2}=0$.

Clearly if we have solutions to (\ref{ve1}), from the third equation
in (\ref{ve1}) we
obtain that $\deg L\leq0$, while from (\ref{ve2}) we have $\deg L\geq0$,
where $\deg L$ is the degree of $L$ with respect to the \kahler\ metric
defined as in  (\ref{degree}).
Since we are interested only in irreducible solutions,
obviously these two situations cannot  occur simultaneously.
The Hodge star operator interchanges these two
cases, and we can thus concentrate on case (i).

Equations (\ref{ve1})  are essentially the equations known as  the
{\em vortex equations}. These are generalisations of the
vortex equations on the Euclidean plane studied by Jaffe and Taubes
\cite{T1,T2,JT}, and have been extensively studied (e.g. in
 \cite{B1,B2,G2,G3} ) for compact \kahler\
manifolds of arbitrary dimension.  In that setting, the equations are the
following:

Let   $(X,\omega)$ be  a compact \kahler\ manifold
of arbitrary dimension, and
let $(L,h)$ be a Hermitian $C^\infty$ line bundle over $X$.
Let $\tau\in\bR$.
The  $\tau$-{\em vortex equations}
\be
\left. \begin{array}{l}
F_A^{0,2}=0\\
\dbar_A\phi=0 \\
\Lambda F_A =\frac{i}{2}(|\phi|^2 -\tau)
\end{array}\right \}, \label{ve}
\ee
are equations for a pair $(A,\phi)$
consisting of a connection on $(L,h)$  and a smooth section of $L$.
The first equation means  that $A$ defines a holomorphic
structure on $L$, while  the second says that $\phi$ must be holomorphic
with respect to this holomorphic structure.

Coming back to  the monopole equations,
we first observe that (\ref{ve1}) can be rewritten  in terms of
$\hA$ only, i.e. not
involving simultaneously $A$ and $\hA$. To do this we recall that
$\hA=(A\otimes a_K)^{1/2}$, where $a_K$ is the metric connection on $K$.
We thus have
\be
F_{\hA}=\frac{1}{2}(F_A + F_{a_K}),  \label{curvature}
\ee
and hence $F_A^{0,2}=0$ is equivalent to $F_{\hA}^{0,2}=0$ since
$F_{a_K}^{0,2}=0$.

 From (\ref{curvature}), and using that $s=-i\Lambda F_{a_K}$ is the scalar
curvature, we obtain that
$$
\Lambda F_A=2\Lambda F_\hA-is,
$$
and hence (\ref{ve1}) is equivalent to
\be
\left. \begin{array}{l}
F_\hA^{0,2}=0\\
\dbar_\hA\phi=0\\
\Lambda F_\hA=\frac{i}{2}(|\phi|^2+ s)
\end{array}
\right\}.\label{ve1'}
\ee

These are the vortex equations on $\hL$, but with the parameter
$\tau$ replaced by minus the scalar curvature. As we will
explain in Section \ref{t-vortices}
the existence proofs for the vortex equations can be
easily modified to give an existence
theorem for the equations obtained by replacing the parameter $\tau$ by
a function $t\in C^\infty (X,\bR)$ in (\ref{ve}).
%%%%%!!!!!!!! CHANGES STARTIN HERE
However, to  compute  the \sw\ invariants one can  slightly perturb
equations (\ref{me})  in such a way that, when $\beta=0$,  equations
(\ref{me}) reduce to the constant function vortex equations, i.e. to
(\ref{ve}) (see e.g. \cite{G5}).
%%%%%%%!!!!ENDING HERE

%%%%%%%%%%%%%%%%%%%%%%%%%%%%%%%%%%%%%%%%%%%%%%%%%%%%%%%%
\section{Non-abelian vortex equations}\label{vortices}
%%%%%%%%%%%%%%%%%%%%%%%%%%%%%%%%%%%%%%%%%%%%%%%%%%%%%%%%
As we have seen in the previous section, the \sw\ monopole
equations can be considered as a four-dimensional Riemannian
generalisation of the vortex equations.
This suggests that we may find  interesting \sw-type equations
by considering the corresponding analogues of
different equations of vortex-type existing in the literature.
With this objective in mind, in this section we shall
review three different non-abelian generalisations of the vortex
equations described above.
The first one consists in studying the vortex equations on
a Hermitian vector bundle of arbitrary rank. The other two
involve two vector bundles, one of which  will  actually be a line
bundle in most cases.

%%%%%%%%%%%%%%%%%%%%%%%%%%%%%%%%%%%%%%%%%%%%
%\subsection{The higher rank vortex equations}
%%%%%%%%%%%%%%%%%%%%%%%%%%%%%%%%%%%%%%%%%%%%
Let $(E,H)$ be a Hermitian vector bundle over a compact
\kahler\ manifold $(X,\omega)$ of complex dimension $n$.
Let $\tau\in\bR$.
 The $\tau$-vortex equations
were generalised to this situation in \cite{B2}. As in the line
bundle case, one studies equations
\be
\left. \begin{array}{l}
F_A^{0,2}=0\\
\dbar_A\phi=0 \\
\Lambda F_A =i(\phi\otimes\phi^\ast -\tau\bI)
\end{array}\right \} \label{nave}
\ee
for a pair $(A,\phi)$ consisting of a unitary connection on $(E,H)$
and a smooth section of $E$.
By $\phi^\ast$ we denote the adjoint of $\phi$ with respect to $H$, and
$\bI\in\End E$ is the identity.
Notice that in the abelian equations (\ref{ve}) there is a $1/2$ in the
third equation, while in (\ref{nave}) there is not. This is not essential
since by applying a constant complex gauge transformation we can introduce
an arbitrary positive constant in front of $i\phi\otimes\phi^\ast $.

These vortex equations appear naturally as the equations satisfied by the
minima of the \ymh\ functional. This is a functional
defined on the product of the space $\cA$ of unitary connections
on $(E,H)$ and the space of smooth sections $\Omega^0(E)$ by
$$
\YMH_\tau(A,\phi)=\|F_A\|^2 +2\|d_A\phi\|^2+
\|\phi\otimes\phi^\ast-\tau\bI\|^2,
$$
where $\|\;\;\|$ denotes the $L^2$-metric.

This is easily seen by  rewriting the \ymh\ functional --- using
the \kahler\ identities --- as

\begin{eqnarray}
\YMH_\tau(A,\phi) &=& 4\|F_A^{0,2}\|^2+4\|\dbar_A\phi\|^2+
\|i\Lambda F_A +\phi\otimes\phi^\ast -\tau\bI\|^2 \nonumber \\
 & &+4\pi\tau\deg E-
\frac{8\pi^2}{(n-2)!} \int_X\ch_2(E)\wedge\omega^{n-2},\nonumber
\end{eqnarray}
where $\deg E$ is the degree of $E$ defined as
\be
\deg E=\int_X c_1(E)\wedge \omega^{n-1}\ , \label{degree}
\ee
and $\ch_2(E)$ is the second Chern character of $E$, which is represented
in terms of the curvature by
$$
\ch_2(E)=-\frac{1}{8\pi^2} \tr(F_A\wedge F_A).
$$
Clearly $\YMH_\tau$ achieves its minimum value
$$
4\pi\tau\deg E-
\frac{8\pi^2}{(n-2)!} \int_X\ch_2(E)\wedge\omega^{n-2}
$$
if and only if $(A,\phi)$ is a solution to equations (\ref{nave})
(see \cite{B1} for details).

As we will explain in Section \ref{t-vortices},  the vortex equations
also have a symplectic interpretation as moment map equations.
The moment map in question is for a symplectic action of $\cG(E)$, i.e.
the unitary gauge group of $E$, on a certain infinite dimensional
symplectic space.
%%%%%%%%%%%%%%%%%%%%%%%%%%%%%%%%%%%%%%%%%
%\subsection{The coupled vortex equations}
%%%%%%%%%%%%%%%%%%%%%%%%%%%%%%%%%%%%%%%%%

A natural generalization of these equations is obtained if we regard the
section $\phi$ in (\ref{nave}) as a morphism from
the trivial line bundle to $E$. One can
replace the trivial line bundle by a vector bundle of arbitrary
rank and study equations for connections on both bundles and
a morphism from one to the other.
These are the {\em coupled vortex equations} introduced in \cite{G4}.

Let $(E,H)$ and $(F,K)$ be smooth Hermitian vector bundles over $X$.
Let $A$ and $B$ be unitary connections on $(E,H)$ and $(F,K)$
resp., and let $\phi\in \Omega^0(\Hom(F,E))$.
Let $\tau$ and $\tau'$ be real parameters.
The coupled vortex equations are
  \be
 \left. \begin{array}{l}
F_A^{0,2}=0\\
F_B^{0,2}=0\\
\dbar_{A, B}\phi=0\\
 i  \Lambda F_A+\phi\phi^\ast=\tau \bI_E\\
i \Lambda F_B-\phi^\ast\phi=\tau'\bI_F
  \end{array}\right \}.\label{gcve}
  \ee

As in the case of the vortex equations described above, equations
(\ref{gcve})
correspond to the minima of a certain \ymh\ functional and are also
moment map equations (cf. \cite{G4}).
The appropriate functional in this case is defined as
$$
\YMH_{\tau,\tau'}(A,B,\phi)=\|F_A\|^2 +\|F_B\|^2 +
2\|d_{A\otimes B}\phi\|^2+
\|\phi\phi^\ast-\tau\bI_E\|^2+\|\phi^\ast\phi-\tau'\bI_F\|^2.
$$

The moment map is now for a symplectic action of $\cG(E)\times\cG(F)$, i.e.
for
 the product of the unitary groups of $E$ and $F$.

In this paper we will be mostly interested in the case in which $F=L$ is
a line bundle. In this situation the equations can be written as
\be
  \left. \begin{array}{l}
F_A^{0,2}=0\\
F_B^{0,2}=0\\
\dbar_{A, B}\phi=0\\
i  \Lambda F_A+\phi\otimes\phi^\ast=\tau \bI_E\\
i \Lambda F_B-|\phi|^2=\tau'
\end{array}\right \}.\label{cve}
\ee
It is clear, from taking the trace of the last two equations
in (\ref{cve}) and integrating, that to solve (\ref{cve})
$\tau$ and $\tau'$ must be related by
\be
\tau\rank E  +\tau'=\deg E +\deg L,\label{parameters}
\ee
hence there is only one free parameter.
%%%%%%%%%%%%%%%%%%%%%%%%%%%%%%%%%%%%%%%%
%\subsection{The framed vortex equations}
%%%%%%%%%%%%%%%%%%%%%%%%%%%%%%%%%%%%%%%%

We  shall consider next the {\em framed vortex equations}. The situation
 is very similar to the previous one in that it also involves  two
vector bundles. In this case, however, the connection on one of the
bundles is fixed. More specifically, let $(E,H)$ and $(F,K)$ be
two Hermitian vector bundles. Let $B$ a fixed Hermitian connection
on $(F,K)$ such that $F_B^{0,2}=0$.

The equations are now for a unitary connection $A$  on $(E,H)$, and
$\phi\in \Omega^0(\Hom(F,E))$.  As explained in \cite{BDGW}, the
appropriate equations are

\be
 \left. \begin{array}{l}
F_A^{0,2}=0\\
\dbar_{A, B}\phi=0\\
 i  \Lambda F_A+\phi\phi^\ast=\tau \bI_E\\
\end{array}\right \}\label{fve}
  \ee
The relation between these equations and the full coupled vortex equations is
perhaps best understood from the symplectic point of view. One sees that the
effect of fixing the data on $F$ is to reduce the symmetry group in the
problem
from $\cG(E)\times\cG(F)$\ to $\cG(E)$. The new equations must thus correspond
to the moment map for the subgroup
$\cG(E)\subset\cG(E)\times\cG(F)$. But the moment maps for $\cG(E)$ and for
$\cG(E)\times\cG(F)$ are related by a projection from the Lie algebra of
$\cG(E)\times\cG(F)$ onto the summand corresponding to the Lie algebra of
$\cG(E)$.  The effect of this projection is to eliminate the last equations in
(\ref{cve}) (cf.
\cite{BDGW} for more details).

The appropriate moduli space problem corresponds to that of studying
morphisms from a vector bundle with a fixed holomorphic structure
to another  vector bundle in which the holomorphic structure is varying.
Such moduli spaces have been studied by Huybrecht and Lehn
\cite{HL1,HL2}, who refer to these objects as {\em framed modules}.

As in  the coupled vortex equations, we will be mostly interested in the
case in which $F=L$ is a line bundle.

All the vortex-type equations that we have considered so far involve
one or two real parameters $\tau$ and $\tau'$. As we have seen in
Section \ref{sw}, the study of the \sw\ monopole equations leads to
abelian vortex
equations in which $\tau$ is replaced by a function. The same will
happen in the generalizations of the monopole equations that we are
about to discuss. This will be analysed in detail in Section
\ref{t-vortices}.

%%%%%%%%%%%%%%%%%%%%%%%%%%%%%%%%%%%%%%%%%%%%%%%%%%%%%%%%%%%
\section{Non-abelian monopole equations}\label{monopoles}
%%%%%%%%%%%%%%%%%%%%%%%%%%%%%%%%%%%%%%%%%%%%%%%%%%%%%%%%%%%
Let us go  back to the set-up of Section \ref{sw} and let  $(X,g)$ be a
compact, oriented, Riemannian,  four-dimensional manifold.
Let $c\in \spinc(X)$ be a fixed $\spinc$-structure, with corresponding
Hermitian line bundle $L$ and bundles of spinors $\slpm$.
Let $(E,H)$ be a Hermitian vector bundle on $X$.
Let $\Psi\in \Gamma (S^+_L\otimes E)$. Using the metrics on $S_L^+$ and
$E$ one has the  antilinear identification
$$
S^+_L\otimes E\lra S^{+*}_L\otimes E^*
$$
$$
\Psi \mapsto \Psi^*.
$$
and hence
$$
\Psi \otimes \Psi^*\in \End (S^+_L\otimes E).
$$
We shall also need the map
$$
\End (S^+_L\otimes E)\lra \End_0(S^+_L)\otimes \End E
$$
$$
\Psi\otimes \Psi^*\mapsto (\Psi\otimes \Psi^*)_0\; ,
$$
as well as the map
$$
\End_0(S^+_L)\otimes \End E\stackrel{\tr}{\lra}\End_0 (S^+_L)
$$
$$
(\Psi\otimes \Psi^*)_0\mapsto \tr(\Psi\otimes \Psi^*)_0
$$
obtained from the trace map
$$
\End_0\,E\lra \End E \stackrel{\tr}{\lra}\bC\lra 0.
$$
The endomorphism  $(\Psi\otimes \Psi^\ast)_0$ should not be confused
with the completely traceless part of
$\Psi\otimes \Psi^\ast$ --- here we are only
removing the trace corresponding to $S_L^+$.

In this paper we shall  consider essentially two different non-abelian
generalizations of the \sw\ equations. While in the first one we will fix
a connection $b$ on $L$
and  study equations  for a pair $(A,\Psi)$, where $A$
is a unitary connection on $(E,H)$ and $\Psi\in \Gamma(S_L^+\otimes E)$,
in the second one we will allow $b$ to vary as well,
and hence our system of equations will be one for triples $(A,b,\Psi)$.

%%%%%%%%%%%%%%%%%%%%%%%%%%%%%%%%%%%%%%%%%%%%%%%%%%%%%%%%%%%%
%\subsection{Monopole equations with fixed connection on $L$}
%%%%%%%%%%%%%%%%%%%%%%%%%%%%%%%%%%%%%%%%%%%%%%%%%%%%%%%%%%%%
The first  non-abelian version of the \sw\ equations is suggested by the
framed vortex equations:
Let $b$ be a {\em fixed} unitary connection on $L$ and  $A$ be a unitary
connection on $E$. Using these two connections and
the Levi--Civita connection one can  consider  the coupled Dirac operator
$$
D_{A,b}:\; \Gamma(\slp\otimes E)\lra \Gamma(\slm\otimes E).
$$
and  study the equations
\be
\left.\begin{array}{l}
D_{A,b}\Psi=0\\
F^+_A=i(\Psi\otimes \Psi^*)_0
\end{array}
\right \}\label{name}
\ee
for the unknowns $A$ and $\Psi\in\Gamma(\slp\otimes E)$.

The Bochner--Weitzenb\"ock formula for $D_{A,b}$ is given by
\be
D_{A,b}^\ast D_{A,b}=
\nabla_{A,b}^\ast\nabla_{A,b}+\frac{s}{4} +c(F_{A,b}), \
\label{bw}
\ee
where
$$
F_{A,b}=F_A +\frac{1}{2} F_b\otimes I_E
$$
and $\nabla_{A,b}$ is the connection on $E\otimes S_L$ determined
by the connections $A$ and $b$.
The term $c(F_{A,b})$ in (\ref{bw}) means Clifford multiplication by
$F_{A,b}$. In fact the action of $F_{A,b}$ on
$\Psi\in \Gamma(\slp\otimes E)$ coincides
with the Clifford multiplication with $F_{A,b}^+$
(see \cite{LaMi} for example).

It is thus  natural to perturb the second equation in (\ref{name})
by adding the constant self-dual 2-form $\frac{1}{2} F_b^+$,
and consider instead equations
\be
\left.\begin{array}{l}
D_{A,b}\Psi=0\\
F^+_{A, b}=i(\Psi\otimes \Psi^*)_0
\end{array}
\right \}. \label{name'}
\ee
These are the equations studied  in \cite{OT1}.

The abelian \sw\ equations (\ref{me}) for a $\spinc$-structure
with line bundle $\tilde L$ can be recovered from (\ref{name})
or (\ref{name'}) by considering the Hermitian bundle
$E=(\tilde L \otimes L^\ast)^{1/2}$.
Notice that this square root exists since
$c_1(\tilde L)\equiv c_1(L)\;\mod\; 2$.

Other non-abelian versions of the \sw\ equations that have been considered
include replacing the $\U(r)$-bundle $(E,H)$ by an $\su(r)$-bundle,
and study  equations  (\ref{name}) in which
$(\Psi\otimes\Psi^\ast)_0$ is replaced by the completely trace-free part
of $\Psi\otimes\Psi^\ast$. These are the equations studied
in \cite{LM} (see also \cite{OT2})
In other versions, like the one considered in \cite{PT}, one fixes the
connection of $\det E\otimes L$ instead of fixing that of $L$.

In all the versions mentioned above one considers spinors coupled to a
bundle associated the fundamental representation of $\U(r)$ or $\su(r)$.
Another direction in which the monopole equations can be generalized is
by considering any compact Lie group $G$ and/or an arbitrary
representation. When the manifold is K\"ahler some of these correspond
to the vortex-type equations described in \cite{G1}. These more general
equations will be dealt with somewhere else.

%%%%%%%%%%%%%%%%%%%%%%%%%%%%%%%%%%%%%%%%%%%%%%%%%%%%%%%%%%%%%%%
%\subsection{Monopole equations with variable connection on $L$}
%%%%%%%%%%%%%%%%%%%%%%%%%%%%%%%%%%%%%%%%%%%%%%%%%%%%%%%%%%%%%%%

We shall consider next the case in which the connection on $L$ is
also varying.
It is clear that we cannot consider the same equations as in the previous
situation,
since we would not obtain  an elliptic complex in the linearization
of the equations --- we need an extra equation.
As the coupled vortex equations (\ref{cve}) suggest, it is natural to
consider the following set of equations for the triple  $(A,b,\Psi)$:
\be
\left.
\begin{array}{l}
 D_{A,b} \Psi =0\\
 F^+_A=i(\Psi \otimes\Psi^*)_0\\
 F^+_b=2i\tr(\Psi\otimes \Psi^*)_0
\end{array}
\right\}.\label{cme}
\ee

%%%%%%%%%%%%%%%%%%%%%%%%%%%%%%%%%%%%%%%%%%
\section{The K\"ahler case}\label{kaehler}
%%%%%%%%%%%%%%%%%%%%%%%%%%%%%%%%%%%%%%%%%%

Let now   $(X,\omega)$ be a compact \kahler\ surface.
Let us fix a $\spinc$-structure $c\in \spinc(X)$, and let $L=L_c$ be the
corresponding Hermitian line bundle.
As mentioned in Section \ref{sw}, the corresponding spinor bundles for the
$\spinc$-structure $c$ are given by
$$
\slp=\hL\oplus \Lambda^{0,2}\otimes \hL\;\;\;\mbox{and}\;\;\;
\slm=\Lambda^{0,1}\otimes\hL,
$$
where
$$
\hL=(K\otimes L)^{1/2}.
$$

Let $(E,H)$ be a Hermitian vector bundle over $X$.
\be
E\otimes\slp=E\otimes\hL\oplus \Lambda^{0,2}\otimes
E\otimes\hL\label{E-spinors}
\ee
and hence
$$
\Gamma(E\otimes\slp)=\Omega^0(E\otimes\hL)\oplus\Omega^{0,2}(E\otimes\hL).
$$
We can write $\Psi$ according to this decomposition as a pair
$\Psi=(\phi,\beta)$.
Let  $b$ and $A$ be  unitary connections on $L$  and  $E$, respectively.
Let $a_K$ be the metric connection on $K$. We shall denote
by $\hb$ the connection on $\hL$ defined by $b$ and $a_K$.
The Dirac operator $D_{A,b}$ can be written in this language as
$$
\dbar_{A,\hb} +\dbar_{A,\hb}^\ast\;:
\Omega^0(E\otimes\hL)\oplus\Omega^{0,2}(E\otimes\hL)\lra
\Omega^{0,1}(E\otimes\hL),
$$
where $\dbar_{A,\hb}$ is the $\dbar$ operator on $E\otimes\hL$
corresponding to the connections $A$ and $\hb$.

%%%%%%%%%%%%%%%%%%%%%%%%%%%%%%%%%%%%
\subsection{Fixed connection on $L$}
%%%%%%%%%%%%%%%%%%%%%%%%%%%%%%%%%%%%
We shall perturb equations (\ref{name'}) by a self-dual 2-form $\alpha$
and consider
\be
\left.\begin{array}{l}
D_{A,b}\Psi=0\\
F^+_{A, b}=i((\Psi\otimes \Psi^*)_0+\alpha \bI)
\end{array}
\right \}. \label{pname}
\ee
We will choose the perturbation to be of \kahler\ type, that is
$\alpha= -f\omega$, where $f$ is a smooth real function.

Similarly to the abelian case, we can write (\ref{pname}) as

\be
\begin{array}{l}
\dbar_{A,\hb}\phi+\dbar^\ast_{A,\hb}\beta=0\\
\Lambda F_{A,b}=i(\phi\otimes \phi^\ast- \Lambda^2
\beta\otimes\beta^\ast-f\bI) \\
F_{A,b}^{2,0}=-\phi\otimes\beta^\ast \\
F_{A,b}^{0,2}=\beta\otimes\phi^\ast\ .
\end{array}\label{kname}
\ee
By  $\Lambda^2$ we denote the operation of contracting twice  with the
\kahler\ form.

As in the abelian case, one can see that the solutions to these equations
are such that either $\beta=0$ or $\phi=0$. More precisely

\begin{prop}\label{decoupling}
Let $\overline{f}=\frac{1}{2\pi}\int_X f$.
The only solutions to  (\ref{kname}) satisfy either

\n (i)\ $\beta=0$,
\be
\begin{array}{l}
F_{A,b}^{0,2}=0\\
\dbar_{A,\hb}\phi=0\\
\Lambda F_{A,b}=i(\phi\otimes\phi^\ast-f\bI)\ ,
\end{array}\label{nave1}
\ee
then $\mu(E)-1/2\deg L\leq\overline{f}$,
or

\n (ii)\ $\phi=0$,
\be
\begin{array}{l}
F_{A,b}^{0,2}=0\\
\dbar_{A,\hb}^\ast\beta=0\\
\Lambda F_{A,b}=-i(\Lambda^2\beta\otimes\beta^\ast+f\bI)
\end{array}\label{nave2}
\ee
and then $\mu(E)-1/2\deg L\geq\overline{f}$.
\end{prop}
\pf
One uses exactly the same method as the one used by Witten \cite{W} in
the abelian case.
Consider the transformation $(A,\phi,\beta)\mapsto (A,-\phi,\beta)$.
Although this is not a symmetry of  equations (\ref{kname}), if
$(A,\phi,\beta)$ is a solution $(A,\phi,-\beta)$ must also be.
This is easily seen by considering the functional
$$
\SW(A,\Psi)=\|F^+_{A, b}-i((\Psi\otimes \Psi^*)_0+\alpha \bI)\|^2
	     +2\|D_{A,b}\Psi\|^2.
$$
Using the Bochner--Weitzenb\"ock formula (\ref{bw}) and the fact that
on  a \kahler\ manifold the decomposition (\ref{E-spinors}) is parallel
with respect to the connection $\nabla_{A,b}$, we have
\begin{eqnarray}
\SW(A,\Psi)& = &\|F_{A,b}^+\|^2+
	    2\|\nabla_{A,b}\phi\|^2+2\|\nabla_{A,b}\beta\|^2+
	    \|i((\Psi\otimes \Psi^*)_0+\alpha \bI)\|^2 \nonumber\\
    & &       + \int_X\frac{s}{2}(|\phi|^2 +|\beta|^2)-
	     2\int_X\langle F_{A,b},i\alpha\bI\rangle.\nonumber
\end{eqnarray}
\remark Notice the analogy between this and the way of rewriting the
\ymh\ functional in Section \ref{vortices} using the \kahler\ identities.
In fact in the \kahler\ case both things are essentially equivalent.
\qed

Of course the only way in which the two type of solutions
can occur simultaneously is if $\mu(E)=1/2\deg L +\overline{f}$, then
$\Psi=0$ and the equations reduce essentially to the \he\ equations.

Since the Hodge operator interchanges the roles of $\phi$ and $\beta$ we
may concentrate in the case $\phi\neq 0$.
We shall write equations
(\ref{nave1}) in a way that we can identify them as the vortex equations
discussed in Section \ref{vortices}. To do this denote by
$$
\hE=E\otimes \hL\;\;\;\mbox{and}\;\;\; \hA=A\otimes \hb,
$$
where recall that $\hb=b\otimes a_K$. We have that
$$
F_\hA=F_{A,b}+\frac{1}{2} F_{a_K},
$$
and hence (\ref{nave1}) is equivalent to
\be
\begin{array}{l}
F_\hA^{0,2}=0\\
\dbar_\hA\phi=0\\
\Lambda F_\hA=i(\phi\otimes\phi^\ast+(s/2-f) \bI)
\end{array}
\label{nave1'}
\ee
where $s=-i\Lambda F_{a_K}$ is the scalar curvature of $(X,\omega)$, and
we have used that $a_K$ is integrable, i.e. $F_{a_K}^{0,2}=0$.
These are indeed the vortex equations (\ref{nave}) on the bundle $\hE$,
with the parameter $\tau$ replaced by the function $t=f-s/2$.

Equations (\ref{nave1}) can also be interpreted as the framed vortex
equations (\ref{fve}). To see this we choose the fixed connection $b$
to be integrable. We then have that $F_{A,b}^{0,2}=F_A^{0,2}$ and
(\ref{nave1}) is equivalent to
 \be
\begin{array}{l}
F_A^{0,2}=0\\
\dbar_{A,\hb^\ast}\phi=0\\
\Lambda F_A=i(\phi\otimes\phi^\ast- t \bI)
\end{array}
\label{nave1''}
\ee
where $t=f+\frac{i}{2}\Lambda F_b$. We thus obtain the framed vortex
equations (\ref{fve}) on $(E,\hL^\ast)$ with fixed connection $\hb^\ast$
on $\hL^\ast$, and parameter $\tau$ replaced by the function $t$.

These two slightly different points of view
in relating (\ref{nave1}) to the vortex equations reflect the close
relation between the usual  vortex equations  and the framed
vortex equations, as we will explain in Section \ref{t-vortices}.
%%%%%%%%%%%%%%%%%%%%%%%%%%%%%%%%%%%%%%%
\subsection{Variable connection on $L$}
%%%%%%%%%%%%%%%%%%%%%%%%%%%%%%%%%%%%%%%
We come now to equations  (\ref{cme}). In the \kahler\ situation
these equations can be  written as
\be
\begin{array}{l}
\dbar_{A,\hb}\phi+\dbar^\ast_{A,\hb}\beta=0\\
\Lambda F_A=i(\phi\otimes \phi^\ast- \Lambda^2
\beta\otimes\beta^\ast) \\
\Lambda F_b=2i(|\phi|^2-|\beta|^2) \\
F_A^{2,0}=-\phi\otimes\beta^\ast\\
F_A^{0,2}=\beta\otimes\phi^\ast\\
F_b^{2,0}=-2\tr(\phi\otimes\beta^\ast)\\
F_b^{0,2}=2\tr(\beta\otimes\phi^\ast) \ . \label{kcme}
\end{array}
\ee
By taking the  third  equation, subtracting twice the trace of the
second equation in (\ref{kcme}), and integrating we obtain that,
in order to have solutions, we need
$$
\deg E = \frac{1}{2}\deg L.
$$
To avoid this restriction we can  perturb, as in the previous case,
the coupled monopole
equations  by adding fixed self-dual forms $\alpha, \gamma$, i.e. by
 considering
\be
\left.
\begin{array}{l}
 D_{A,b} \Psi =0\\
 F^+_A=i(\Psi \otimes\Psi^*)_0+i\alpha \bI_E\\
 F^+_b=2i\tr(\Psi\otimes \Psi^*)_0 +i\gamma \label{pcve}
\end{array}
\right\}.
\ee
In the \kahler\ case we shall choose
$$
\alpha=-f\omega \;\;\;\;\;\;\;\gamma=2f'\omega,
$$
where $f$ and $f'$ are  smooth real functions on $X$.

With this choice of perturbation  the second  and
third equations in (\ref{kcme}) become respectively
$$
i \Lambda F_A + \phi \otimes \phi^\ast -
\Lambda^2 \beta \otimes\beta^* = f \bI_E
$$
$$
i \Lambda F_b + 2(|\phi|^2 - |\beta|^2) =- 2f'.
$$
A necessary condition for existence of solutions is now
\be
\deg E - \frac{1}{2}\deg L =  \int_X (rf + f'), \label{ff'}
\ee
where $r = \rank E$.

We shall study the more general monopole equations (\ref{pcve}).

\begin{prop}
Let  $f,f'\in C^\infty(X,\bR)$ be  related by (\ref{ff'}) and
denote
$$
\overline{f}=\frac{1}{2\pi}\int_X f \;\;\;\;\mbox{and}\;\;\;\;
\overline{f'}=\frac{1}{2\pi}\int_X f'.
$$
The only solutions to (\ref{kcme}) are such that either

\n (i)\ $\beta=0$,
\be
 \begin{array}{l}
F_A^{0,2}=0\\
F_b^{0,2}=0\\
\dbar_{A,\hb}\phi=0\\
i\Lambda F_A+\phi\otimes\phi^\ast=f\bI_E\\
i\Lambda F_b+2|\phi|^2=-2f'\ ,
\end{array}
\label{cve1}
\ee
then $\mu(E)\leq \overline{f}$, $\deg L\geq 2\overline{f'}$
or

\n (ii)\ $\phi=0$,
\be
\begin{array}{l}
F_A^{0,2}=0\\
F_b^{0,2}=0\\
\dbar_{A,\hb}^\ast\beta=0\\
i\Lambda F_A-\Lambda^2\beta\otimes\beta^\ast=f  \bI_E\\
i\Lambda F_b-2|\beta|^2=-2f'
\end{array}\label{cve2}
\ee
and then $\mu(E)\geq \overline{f}$, $\deg L\leq 2\overline{f'}$.
\end{prop}
\pf
This is proved similarly to Proposition \ref{decoupling}, by considering
the functional
$$
\SW(A,b,\Psi)=
\| F^+_A-i(\Psi \otimes\Psi^*)_0-i\alpha \bI_E\|^2+
 2\|F^+_b-2i\tr(\Psi\otimes \Psi^*)_0 -i\gamma\|^2+
 2\| D_{A,b} \Psi\|^2,
$$
and observing that
$$
\langle F_A^+,i(\Psi \otimes\Psi^*)_0\rangle +\frac{1}{2}
\langle F_b^+,2i\tr(\Psi \otimes\Psi^*)_0\rangle =
\langle F_{A,b}^+,i(\Psi \otimes\Psi^*)_0\rangle.
$$
\qed

Again, we can focus on case (i), since by means of the Hodge operator
we can interchange the roles of $\phi$ and $\beta$.

The system of equations (\ref{cve1}) is equivalent to
\be
\begin{array}{l}
F_A^{0,2}=0\\
F_{\hb^\ast}^{0,2}=0\\
\dbar_{A,\hb^\ast}\phi=0\\
i\Lambda F_A+\phi\otimes\phi^\ast=f\bI_E\\
i\Lambda F_{\hb^\ast}-|\phi|^2=(f'-\frac{s}{2})
\end{array}
\label{cve1'}
\ee
where $\hb ^\ast$ is the dual connection to $\hb$ on $\hL^\ast$, which
satisfies that $F_{\hb^\ast}=-F_\hb$.

We can now identify (\ref{cve1'}) as the coupled vortex equations
(\ref{cve}) on $(E,\hL^\ast)$, with the  parameters $\tau$ and $\tau'$
replaced by functions.
The existence of solutions to these equations as well as the description
of the moduli space of all solutions will be dealt with in the next
section.
%%%%%%%%%%%%%%%%%%%%%%%%%%%%%%%%%%%%%%%%%%%%%%%%%%%%%%%%
\section{Back to the vortex equations}\label{t-vortices}
%%%%%%%%%%%%%%%%%%%%%%%%%%%%%%%%%%%%%%%%%%%%%%%%%%%%%%%%
We study now the existence of solutions to the monopole
equations in the \kahler\ case --- equations (\ref{nave1'}),
(\ref{nave1''}) and
(\ref{cve1'}) --- or equivalently the vortex equations in which
the parameters have been replaced by functions.

%%%%%%%%%%%%%%%%%%%%%%%%%%%%%%%%%%%%%%%%%%%%%%%%%
\subsection{The  $t$-vortex equations}
%%%%%%%%%%%%%%%%%%%%%%%%%%%%%%%%%%%%%%%%%%%%%%%%%
In this section  we will examine what happens if the parameter
$\tau$, which appears on the $\tau$-vortex equations (\ref{nave})
is permitted to be a non-constant smooth function, say $t$.
We will refer to the resulting equations as the $t$-vortex equations.
The main result, namely that existence of solutions is
governed entirely by the average value of $t$, has already
been observed by Okonek and Teleman. Here we give a somewhat
different proof than that in \cite{OT1}. We also discuss some
interpretations and implications of the result.

Let  $E\ra X$\ be a rank $r$, smooth complex bundle over a
closed \kahler\ manifold $(X,\omega)$.
It is  convenient to look at equations (\ref{nave}) from a
different although equivalent point of view.
Instead of fixing a Hermitian metric $H$ and solving
for $(A,\phi)$ satisfying (\ref{nave}) we fix a
$\dbar$-operator on $E$, $\dbar_E$ say, and a section
$\phi\in H^0(X,\cE)$, where
$\cE$ is the holomorphic bundle determined by $\dbar_E$,
and solve for a metric $H$ satisfying
\be
i\Lambda F_H +\phi\otimes\phi^{*_H}=\tau\bI,\label{mve}
\ee
where $F_H$\ is the curvature of the metric connection
 determined by $\dbar_E$\ and $H$. It will be important to explicitly
write $\phi^{*_H}$\ to denote  the adjoint  of $\phi$ with respect to the
metric $H$.

Equation (\ref{mve}) can be regarded as the defining condition
for a special metric on the holomorphic pair $(\cE,\phi)$.
In \cite{B2} the first author showed that there is a \hk\ correspondence
between the existence of such metrics and a stability condition  for
holomorphic pairs. The appropriate notion of stability is as follows.
%%%%%%%%%%%%%%%%%%%%%%%%%%%%%%%%%%%%%%

\begin{definition} Define the degree of any coherent sheaf
$\cE'\subset\cE$ to be
$$
\deg\cE'=\int_X c_1(\cE')\wedge \omega^{n-1},
$$
and define the slope of $\cE'$ by
$$
\mu(\cE')=\frac{\deg\cE'}{\rank \cE'}.
$$
Fix $\tau\in\bR$. The holomorphic pair $(\cE,\phi)$ is called  \ts\
if

(1)\  $\mu(\cE')<\tau\;\;
\mbox{for every coherent subsheaf}\;\;\;\cE'\subset\cE$, and

(2)\ $\mu(\cE/\cE'')>\tau\;\;
\mbox{for every non-trivial coherent subsheaf} \;\cE''\; \mbox{such that}
\; \phi\in H^0(X,\cE'')$.
\end{definition}

The \hk\ correspondence is expressed by the following two propositions.

\begin{prop}[\cite{B2}]\label{easy}
Fix $\tau\in \bR$, let $(\cE,\phi)$\ be a  holomorphic pair, and suppose
that there exists a metric, H, satisfying the $\tau$-vortex equation
(\ref{mve}).   Then
either

(1)\ the holomorphic pair $(\cE,\phi)$\ is $\tau$-stable, or

(2)\ the bundle $\cE$\ splits holomorphically as $\cE=\cE_{\phi}\oplus
\cE_{ps}$\ with $\phi\in H^0(X,\cE_{\phi})$, and such that
the holomorphic pair $(\cE_{\phi},\phi)$\ is $\tau$-stable, and
$\cE_{ps}$\ is polystable with slope equal to $\tau$.
\end{prop}

\begin{prop}[\cite{B2}]\label{hard}
Fix $\tau\in \bR$, let $(\cE,\phi)$\ be a $\tau$-stable holomorphic
pair.
Then there is a unique smooth metric, H, on $E$\ such that the
$\tau$-vortex equation (\ref{mve})\ is satisfied.
\end{prop}

\remark We are assuming that the \kahler\ metric is normalized so that
$\vol(X)=2\pi$. Otherwise we need to introduce the factor
$\frac{2\pi}{\vol(X)}$ in the right hand side of (\ref{mve}).

%%%%%%%%%%%%%%%%%%%%%%%%%%%%%%%%%%%%%
%\subheading{\hfil \S 2. Results}
%%%%%%%%%%%%%%%%%%%%%%%%%%%%%%%%%%%%%

Suppose now that we replace $\tau$\ in equation (\ref{mve}) by a
smooth real valued function $t\in C^{\infty}(X,\bR)$ and study
\be
i\Lambda F_H +\phi\otimes\phi^{*_H}=t\bI.\label{tmve}
\ee
The question we wish to address is: how does this affect
the Hitchin-Kobayashi correspondence?

One direction is clear: replacing the constant $\tau$\ by the smooth function
$t$\
has absolutely no effect on the proof of Proposition \ref{easy}
(Theorem 2.1.6 in \cite{B2}). The same proof thus yields

\begin{thm}[cf. also Theorem 3.3, \cite{OT1}] Fix a smooth function
$t\in C^{\infty}(X,\bR)$.
Let $\overline{t}=\frac{1}{2\pi}\int t$.
If a holomorphic pair $(\cE,\phi)$\
supports a solution to the $t$-vortex equation (\ref{tmve}), then either

(1)\ the holomorphic pair $(\cE,\phi)$\ is $\overline t$-stable, or

(2)\  the bundle $\cE$\ splits holomorphically as
$\cE=\cE_{\phi}\oplus\cE_{ps}$\ with $\phi\in H^0(X,\cE_{\phi})$, and
such that
the holomorphic pair $(\cE_{\phi},\phi)$\ is $\overline t$-stable, and
$\cE_{ps}$\ is polystable with slope equal to $\overline t$.
\end{thm}

We now consider the converse result, i.e. the analogue of
Proposition \ref{hard} (Theorem 3.1.1 in \cite{B2}).   As shown in
\cite{OT1}, the $t$-vortex equation
can be reformulated as an equation with a constant right hand side.
Let $\tau=\overline{t}$. Since
$\int_X(\tau-t)=0$, we can find a smooth function $u$\ such
that $\Delta(u)=\tau-t$.  Thus  (\ref{tmve}) is equivalent to
\be
i\Lambda (F_H+\Delta(u)\bI)
  +\phi\otimes\phi^{*_H}=\tau\bI\ .\label{delta-ve}
\ee

If we define a new metric $K=He^u$,

then (\ref{delta-ve}) becomes
\be
i\Lambda F_K
  +e^{-u}\phi\otimes\phi^{*_K}=\tau\bI\ .\label{uveOT}
\ee

This is {\it almost} the $\tau$-vortex equation, the only difference being the
prefactor $e^{-u}$\ in the second term.  In the analysis of this situation by
Okonek and Teleman, they indicate how the proof (of Theorem 3.1.1) in [B2] can
be modified to accomodate this new wrinkle. (The proof in [B2] employs a
modified Donaldson functional on the space of Hermitian metrics on the bundle
$\cal E$. In [OT1], the authors generalize the functional so that it
accomodates the extra factor of $e^{-u}$, and argue that this has little
effect
on the proof.)

It is interesting to observe that the same result can be achieved without {\it
any modification at all} of the proof in [B2], if one enlarges the catgory in
which the proof is applied. This can be seen as follows. If we set

\be
\phi_u=e^{-u/2}\phi\ ,\label{phiu}
\ee
then (\ref{uveOT}) becomes
\be
i\Lambda F_K
  +\phi_u\otimes\phi_u^{*_K}=\tau\bI\ .\label{uve}
\ee

We thus see that
\begin{lemma} \label{t-tau}
Let $u$\ be given by $\Delta(u)=\tau-t$.
The pair $(\cE,\phi)$\ admits a metric $H$ satisfying the
$t$-vortex equation if and only if
the pair $(\cE,\phi_u)$\ admits a metric $K$ satisfying
the $\tau$-vortex equation. The metrics $H$\ and $K$\ are
related by  $K=He^u$.
\end{lemma}

It is important to notice that, unless $u$\ is a constant
function, $\phi_u$\ is {\it not} \ a holomorphic section of
$\cE$. Indeed $\dbar_E\phi_u=-\frac{1}{2}\dbar(u)\phi_u$.
Our key observation is that in Theorem 3.1.1 in \cite{B2}, the
holomorphicity of $\phi$\ is not required either in the
statement or in the proof of the theorem. The proof, and
thus the result remains unchanged if the holomorphic section
$\phi$\ is replaced by a smooth section $\phi_u$\ related to
$\phi$\ by $\phi_u=e^{-u/2}\phi$. The basic reason can be
traced
back to the following simple fact:

\begin{lemma}\label{u-simple}
Let $\phi$, $u$, and $\phi_u$\ be as above.

(1)\ Let $\cE'\subset\cE$\ be any holomorphic
subbundle of $\cE$. Then $\phi$\ is a section of $\cE'$\ if and only if
$\phi_u$\ is a section of $\cE'$.

(2)\ Let $s$\ be any smooth endomorphism of $\cE$. Then
$s\phi=0$\ if and only if $s\phi_u=0$.
\end{lemma}

 Notice, for instance, that if we define
$$
\begin{array}{ll}
H^0_u(X,\cE)&=e^{-u/2}H^0(X,\cE)\\
&=\{e^{-u/2}\phi\in\Omega^0(X,E)\ |\ \phi\in H^0(X,\cE)\}\ ,
\end{array}
$$
then the definition of $\tau$-stability can be applied to any pair
$(\cE,\phi_u)$\ where $\phi_u\in H^0_u(X,\cE)$.

\begin{definition} A pair $(\cE,\phi_u)$, where
$\phi$\ is in
$H^0_u(X,\cE)$, will be called a  $u$-{\em holomorphic
pair}.
\end{definition}

In view of Lemma \ref{u-simple}, it follows that the $u$-holomorphic
pair $(\cE,\phi_u)$\ is $\tau$-stable if and only if the
holomorphic pair $(\cE,\phi)$\ is $\tau$-stable (where $\phi$\
and $\phi_u$\ are related by (\ref{phiu})).  Furthermore, without any
alteration whatsoever, the proof of Theorem 3.1.1
in \cite{B2} can be
applied to a $u$-holomorphic pair to prove:

\begin{prop}\label{upair-tvor}
Fix $u\in C^{\infty}(X,\bR)$\ and
$\tau\in\bR$. Let $(\cE,\phi_u)$\ be a $\tau$-stable $u$-holomorphic pair.
Then $E$\ admits a unique smooth metric, say $K$, such that the
$\tau$-vortex equation is satisfied, i.e. such that
$$
i\Lambda F_K +\phi_u\otimes\phi_u^{*_K}=\tau\bI\ .
$$
\end{prop}

Taken together,  Lemma \ref{t-tau} and Proposition \ref{upair-tvor} thus
prove

\begin{thm} Fix $\tau=\overline t$\ and suppose that
$(\cE,\phi)$\ is a $\tau$-stable pair. Then $E$\ supports
a metric satisfying the $t$-vortex equation.
\end{thm}

%%%%%%%%%%%%%%%%%%%%%%%%%%%%%%%%%%%%%%%%%%
%\subheading {\hfil \S 3. Interpretations}
%%%%%%%%%%%%%%%%%%%%%%%%%%%%%%%%%%%%%%%%%%

The above results describe the sense in which the vortex equation is
{\em insensitive} to the precise form of the parameter $t$.  This can be
made precise by considering the  moduli spaces. Let $\cC$\  be the
space of
holomorphic structures (or, equivalently, integrable $\dbar_E$-operators)
on $E$, and let
$$
\cH=\{(\dbar_E,\phi)\in\cC\times\Omega^0(X,E)\ |\
\dbar_E\phi=0\}\
$$
be the configuration space of holomorphic pairs on $E$.
Let $\cV_t\subset\cH$\ be the set of  $t$-{\em vortices}, i.e.
$$
\cV_t=\{(\dbar_E,\phi)\in\cH\ |\ \mbox{there is a metric satisfying
the $t$-vortex equation}\}\ .
$$
The above results can then be summarized by saying that

(1)\ $\cV_t=\cV_{\tau} $\ for all functions $t$\ such that
$\overline t=\tau$, and

(2)\ for generic values of $\tau$, we can identify
$\cV_t=\cH_{\tau}$\ where $\cH_{\tau}$\ denotes the set of
\ts\ holomorphic pairs. In fact, as complex spaces,
$\cV_t/\cG^c=\cV_{\tau}/\cG^c=\cB_{\tau}$, where $\cB_{\tau}$\
is the moduli space of $\tau$-stable holomorphic pairs ---
which has the structure of a variety
(cf. \cite{Be,BD1,BD2,G4,HL1,HL2,Th}).

\remark In the case in which $E=L$ is a line bundle the $\tau$-stability
condition reduces to\
$\deg L<\tau$\ , and the moduli space of \ts\ pairs is nothing else but
 the space of {\em non-negative divisors} supported by $L$, where by a
non-negative divisor we mean either an effective divisor or
the zero divisor.

Nevertheless, the function $t$\ does carry some information. For example,
the metrics which satisfy the $t$-vortex equation (for fixed $\dbar_E$\
and $\phi$) do depend on $t$. The following observations shed some light
on the role  played by $t$.
%%%%%%%%%%%%%%%%%%%%%%%%%%%%%%%%%%%%%

As described in \cite{G4}, \cite{BDGW}, the holomorphic pair $(\cE,\phi)$\
can be identified with the {\em holomorphic  triple}
$(\cE,\cO,\phi)$, where $\cO$ is the structure sheaf and $\phi$ is
a morphism $\cO\ra\cE$. From this point of view, the
natural equations to consider  are the framed vortex equations
(\ref{fve}).
Coming back to the set-up of Section  \ref{vortices}, we want to study
equations (\ref{fve}) for a vector bundle $E$ of arbitrary rank and
$F=L_0$, the trivial line bundle.
As for the usual vortex equations, we will look at (\ref{fve}) as
equations for a metric on $E$.
To do this we fix the  holomorphic structure
$\dbar_{L_0}$ on $L_0$ to be the trivial one  i.e.
$(L_0,\dbar_{L_0})=\cO$,
and consider a holomorphic  structure $\dbar_E$ on $E$. Then we take
$\phi:\cO\ra\cE$ to be a holomorphic morphism, where $\cE=(E,\dbar_E)$.
In contrast with the coupled vortex equations (\ref{cve}) that we will
analyse later, here we need to fix a metric $h$ on $L_0$.
Then solving (\ref{fve})
is equivalent to solving for a metric $H$ on $E$ satisfying
\be
i\Lambda F_H+\phi\otimes\phi^\ast=\tau\bI\ .\label{mfve}
\ee
It is important to notice that now  $\phi^\ast$ is the adjoint of
$\phi$ with respect to  both metrics $H$ and $h$.
The identification between  $(\cE,\phi)$\ and
$(\cE,\cO,\phi)$\ requires a choice of trivializing frame for $\cO$,
say $f$.
If  $h(f,f)=e^u$, then (\ref{mfve}) becomes
$$
i\Lambda F_H +e^{-u}\phi\otimes\phi^{*_H}=\tau \bI\ .
$$
Thus we recover the usual $\tau$-vortex equation when we select the
metric on $L_0$\ for which the holomorphic frame  of $\cO$\
is also a unitary frame. For other choices of $h$\ we see that we get
essentially equation (\ref{uve}),
or equivalently, the $t$-vortex equation.

{\em  From this point of view, the
function $t$\ is determined by the metric on $L_0$.}

The impact of non-constant $t$\ can also be understood from the symplectic
point of view.  If we fix a metric, say $H$, on $E$, the induced inner
products
on $\cC$\ and on $\Omega^0(X,E)$\ can be combined to give a symplectic
structure on the configuration space $\cH$. Denoting the symplectic
forms
on  $\cC$\ and $\Omega^0(X,E)$\ by $\omega_{H,\cC}$\ and $\omega_{H,0}$\
respectively, we take
$$
\omega_{H,H}=\omega_{H,\cC}+\omega_{H,0}\
$$
as the symplectic form on $\cH$.  Let $\cG_H$\ be the unitary gauge
group of $E$\ determined by $H$, and let $\lie \cG_H$\ be its Lie algebra. A
moment map
$\mu:\cH\ra\lie\cG_H^*$\ for the action of $\cG_H$\ on the symplectic
space $(\cH,\omega_{H,H})$\ is given by
$$
\mu_{H,H}(\dbar_E,\phi)=\Lambda F_{\dbar_E,H} -i\phi\otimes\phi^{*_H}\ ,
$$
where we have written $F_{\dbar_E,H}$ instead of $F_H$ to emphasize that
$F_H$ depends also on the holomorphic structure on $E$.
The $\tau$-vortex equation is thus equivalent to the condition
$\mu_{H,H}(\dbar_E,\phi)=-i\tau\bI$, and we get an identification of
moduli spaces:
$$
\cB_{\tau}=\cV_{\tau}/\cG^c=\mu_{H,H}^{-1}(-i\tau\bI)/\cG_H\ .
$$

Replacing $\tau$\ by the non-constant function $t$\ has no effect on the
identification $\cV_{t}/\cG^c=\mu_{H,H}^{-1}(-it\bI)/\cG_H$:  The
element $it\bI$\
is still a central element in $\lie\cG_H^*$, so the symplectic quotient
at this level is well defined.
(The problem comes in proving that $\cV_t/\cG^c=\cB_{\tau}$.)
An alternative point of view makes use of the equivalence between
the equations (\ref{tmve}) and (\ref{uve}). We define
$$
\mu_{H,K}(\dbar_E,\phi)=\Lambda F_{\dbar_E,H} -i\phi\otimes\phi^{*_K}\ ,
$$
where $H$\ and $K$\ are metrics on $E$. By the above results, $K=He^u$\
where $\Delta(u)=\tau-t$,  then
$$
\mu_{H,H}^{-1}(-it\bI)=
\mu_{H,K}^{-1}(-i\tau\bI)\ .
$$
The point is that $\mu_{H,K}$\ is also a moment map for the action of
$\cG_H$. It arises when the symplectic structure on $\cH$\ is taken to
be
$$
\omega_{H,K}=\omega_{H,\cC}+\omega_{K,0}\ .
$$

{\em From this point of view, the function $t$\ arises from a deformation
of  the symplectic structure on $\cH$.}

%%%%%%%%%%%%%%%%%%%%%%%%%%%%%%%%%%%%%%%%%%%%
\subsection{The coupled vortex equations}
%%%%%%%%%%%%%%%%%%%%%%%%%%%%%%%%%%%%%%%%%%%%

Let us consider the set-up in Section \ref{vortices} for the \cves\
(\ref{cve}).
As in the previous situation, we want to look at  (\ref{cve})
as equations for metrics. In order to do this let us fix  holomorphic
structures $\dbar_E$ and $\dbar_L$ on $E$ and $L$ respectively. Denote by
$\cE$ and $\cL$ the corresponding holomorphic vector bundles. Let
$\phi\in H^0(\cE\otimes\cL^\ast)$. Equations (\ref{cve}) are then
equivalent to solving
\be
\left. \begin{array}{l}
i  \Lambda F_H+\phi\otimes\phi^\ast=\tau \bI_E\\
i \Lambda F_K-|\phi|^2=\tau'
\end{array}\right \}.\label{mcve}
\ee
for metrics $H$ and $K$ on $\cE$ and $\cL$ respectively.

A \hk\ correspondence was proved in \cite{G4}.
The appropriate notion of stability for $(\cE,\cL,\phi)$
can be expressed in terms of the stability of a pair, namely
\begin{definition}\label{st}
The holomorphic triple $(\cE,\cL,\phi)$ is said to be $\tau$-stable if
the holomorphic pair $(\cE\otimes\cL^\ast,\phi)$ is $(\tau-\deg L)$-stable.
\end{definition}

\begin{thm}[\cite{G4}]\label{existence-cve}
Let $\tau$ and $\tau'$ be  real numbers satisfying
(\ref{parameters}).
Let $(\cE,\cL,\phi)$ a holomorphic triple. Suppose that there exist metrics
$H$ and $K$ satisfying (\ref{mcve}), then either
$(\cE,\cL,\phi)$ is \ts\ or
the bundle $\cE$ splits holomorphically as $\cE_\phi\oplus\cE_{ps}$
with $\phi\in H^0(X,\cE_{\phi}\otimes\cL^\ast)$, and such that
$(\cE_\phi,\cL,\phi)$ is \ts\ and $\cE_{ps}$ is polystable with
slope equal to $\tau$.

Conversely, let $(\cE,\cL,\phi)$ be a \ts\ triple then there are unique
smooth metrics $H$ and $K$ satisfying the coupled vortex equations
(\ref{mcve}).
\end{thm}

Suppose now that we replace $\tau$ and $\tau'$ in (\ref{mcve}) by
smooth functions $t, t'\in C^\infty(X,\bR)$, i.e. we consider
\be
\left. \begin{array}{l}
i  \Lambda F_H+\phi\otimes\phi^\ast=t \bI_E\\
i \Lambda F_K-|\phi|^2=t'
\end{array}\right \}.\label{tmcve}
\ee

The first thing that we observe is that in order to have solutions
$t$ and $t'$ must satisfy
\be
\int_X (r t+t')=\deg E +\deg L\label{t-t'}
\ee
where $r=\rank E$.

Next, let us recall the main ideas in the proof of Theorem
\ref{existence-cve}:
The basic fact is that the coupled vortex
equations (\ref{mcve}) are
a dimensional reduction of the \he\ equation
for a metric on a certain vector bundle over $\xp$. This bundle $\cF$,
canonically associated to the holomorphic triple $(\cE,\cL,\phi)$,
is an extension on $\xp$ of the form
\be
0\lra p^\ast\cE\lra \cF\lra \ps \cL \otimes q^\ast\cO(2)\lra 0,
\label{bigbun}
\ee
where $p$ and $q$ are the projections from $\xp$
to $X$ and $\bP^1$  respectively.
This is simply because
$H^1(\xp, p^\ast (\cE\otimes\cL^\ast)\otimes q^\ast\cO(2))
\cong  H^0(X,\cE\otimes\cL^\ast)$

Let  $SU(2)$ act  on $X\times \bP^1$,
trivially  on $X$, and in the standard way  on ${\bP}^1 \cong
SU(2)/U(1)$. This action can be lifted to an action on $\cF$,
 trivial  on $p^\ast{\cal E}$ and $\ps\cL$, and  standard  on
$q^\ast\cO(2)$. The bundle $\cF$ is in this way an $SU(2)$-equivariant
holomorphic vector bundle.

Let $\tau$ and $\tau'$ be related by (\ref{parameters}) and let
\be
\sigma=\frac{4\pi}{\tau-\tau'}
\ee
be positive.
Consider  the $SU(2)$-invariant K\"{a}hler
metric on $X\times \bP^1$ whose K\"{a}hler form is
$$
\omega_\sigma =p^\ast\omega +\sigma  q^\ast\omega_{\bP^1},
$$
where $\omega$ is the K\"{a}hler form on $X$ (normalized such that
$\vol(X)=2\pi$), and $\omega_{\bP^1}$ is
the {\em Fubini-Study} metric with volume 1.

Theorem \ref{existence-cve} is then a consequence of the following two
propositions
and the \hk\ correspondence proved by Donaldson \cite{D1,D2}, and
Uhlenbeck and Yau \cite{UY}).

\begin{prop}[\cite{G4}]\label{dr}
The triple $(\cE,\cL,\phi)$ admits a solution to
(\ref{mcve}) if and only if the vector bundle $\cF$ in (\ref{bigbun})
has a ($\su(2)$-invariant) \he\ metric with respect to $\omega_\sigma$.
\end{prop}
\begin{prop}[\cite{G4}]\label{s-ts}
Suppose that $\cE$ is not isomorphic to $\cL$. Then the triple
$(\cE,\cL,\phi)$
is \ts\ if and only if $\cF$ is stable with respect to $\omega_\sigma$.
If $\cE\cong\cL$, then
$\cF\cong \ps \cL\otimes\qs\cO(1)\oplus\ps \cL\otimes\qs\cO(1)$.
\end{prop}

Suppose first that the functions $t$ and $t'$, in addition to satisfying
(\ref{t-t'}), verify that there  is a positive constant $\sigma$ so that
\be
t-t'=\frac{4\pi}{\sigma}.\label{strong-t-t'}
\ee
The proof of  Proposition \ref{dr} then yields
\begin{prop} \label{weakdr}
The triple $(\cE,\cL,\phi)$ admits a solution to
(\ref{tmcve}) if and only if the vector bundle $\cF$ in (\ref{bigbun})
has a ($\su(2)$-invariant) metric satisfying the weak \he\
equation
\be
i\Lambda_\sigma F_\bH= t\ \bI
\ee
with respect to the \kahler\ form
$\omega_\sigma=\ps \omega_X\oplus \sigma\qs\omega_{\bP^1}$.

\end{prop}
But  the existence of a weak \he\ metric is in fact  equivalent
to the existence of a \he\ metric --- as one can see simply by applying a
conformal change to the metric --- and hence equivalent to the
stability of the bundle. We can then combine again Propositions
\ref{weakdr} and \ref{s-ts} to prove the following.

\begin{thm}\label{existence-tcve}
Fix smooth functions $t,t'\in C^\infty(X,\bR)$ satisfying (\ref{t-t'})
and (\ref{strong-t-t'}). Let  $\overline{t}=\frac{1}{2\pi}\int_X t$ and
$\overline{t'}=\frac{1}{2\pi}\int_X t'$.
Let $(\cE,\cL,\phi)$ a holomorphic triple. Suppose that there exist
metrics $H$ and $K$ satisfying (\ref{tmcve}), then either
$(\cE,\cL,\phi)$ is $\overline{t}$-stable or
the bundle $\cE$ splits holomorphically as $\cE_\phi\oplus\cE_{ps}$
with $\phi\in H^0(X,\cE_{\phi}\otimes\cL^\ast)$, and such that
$(\cE_\phi,\cL,\phi)$ is $\overline{t}$-stable
and $\cE_{ps}$ is polystable with slope equal to $\overline{t}$.

Conversely, let $(\cE,\cL,\phi)$ be a $\overline{t}$-stable triple then
there are unique smooth metrics $H$ and $K$ satisfying
equations(\ref{tmcve}).
\end{thm}

We will  show now that the general coupled vortex equations (\ref{tmcve})
with $t$ and $t'$ satisfying simply (\ref{t-t'})
are also a  dimensional reduction, but in this case of  a metric on $\cF$
satisfying a certain deformation  of the \he\ condition. We set, as above,
$$
\sigma=\frac{4\pi}{\overline{t}-\overline{t'}}
$$
where  $\overline{t}$ and
$\overline{t'}$ denote the average values of $t$
and $t'$ respectively.  Again the proof of  Proposition \ref{dr}
yields

\begin{prop} \label{newdr}
The triple $(\cE,\cL,\phi)$ admits a solution to
(\ref{tmcve}) if and only if the vector bundle $\cF$ in (\ref{bigbun})
has a ($\su(2)$-invariant) metric satisfying the deformed \he\
equation
\be
i\Lambda_\sigma F_\bH= \overline{t}\ \bI +
\left(\begin{array}{cc}(t-\overline{t}) \bI_1 & 0\\
	     0&(t'-\overline{t'})\bI_2 \end{array}\right),
\label{bigbun-eqn}
\ee
with respect to the \kahler\ form
$$
\omega_\sigma=\ps \omega_X\oplus \sigma\qs\omega_{\bP^1}.
$$
\end{prop}

The deformed \he\ equation in this Proposition is similar to the kind
studied in \cite{BG2}.  In \cite{BG2} we considered an extension of
holomorphic
vector bundles over a compact \kahler\ manifold
\be
\extn\label{extn}
\ee
and studied metrics $H$ on $\cE$ satisfying the equation
\be
i\Lambda F_H=\left(\begin{array}{cc}\tau_1 \bI_1 & 0\\
	     0&\tau_2 \bI_2 \end{array}\right),\label{ext-eqn}
\ee
where $\tau_1$ and $\tau_2$ are real numbers, related by
$$
\tau_1 r_1 +\tau_2 r_2=\deg \cE.
$$
The reason we can write an equation like (\ref{ext-eqn}) is that the
metric $H$ gives a  $C^\infty$ splitting of (\ref{extn}).
We proved an existence theorem for metrics  satisfying
(\ref{ext-eqn}) in terms of a notion of stability for the extension
depending on the parameter $\alpha=\tau_1-\tau_2$.
To define this stability condition consider any coherent subsheaf
$\cE'\subset\cE$ and write it as a subextension
$$
\sextn.
$$
Define the $\alpha$-slope of $\cE'$ as
$$
\mu_\alpha(\cE')=\mu(\cE')+\alpha\frac{\rank\cE_2'}{\rank \cE'}.
$$
Then we say that (\ref{extn}) is $\alpha$-stable if and only if for every
non-trivial subsheaf $\cE'\subset\cE$
$$
\mu_\alpha(\cE')<\mu_\alpha(\cE).
$$
We proved
\begin{thm}[\cite{BG2}]Let $\alpha=\tau_1-\tau_2\leq 0$ and
suppose that (\ref{extn}) is indecomposable (as an extension), then
$\cE$ admits a metric satisfying (\ref{ext-eqn}) if and only
(\ref{extn}) is $\alpha$-stable.
\end{thm}
\remark If $\alpha=0$ (\ref{ext-eqn}) reduces to the \he\ equation
and the stability condition is the usual stability of the bundle
$\cE$.

The deformed \he\ equation in Proposition \ref{newdr} differs from
(\ref{ext-eqn}) only in that the constants $\tau_1$ and $\tau_2$
have been replaced by smooth functions, $t_1$ and $t_2$,
satisfying
\be
\int(r_1 t_1 +r_2 t_2)=\deg \cE.\label{t1-t2}
\ee

By the same methods used in \cite{BG2} one can readily show one direction of
the
Hitchin-Kobayashi correspondence, namely

\begin{thm}\label{alpha-hk}Let\  $t_1$ and $t_2$ be  smooth real functions
satisfying (\ref{t1-t2}) and such that $\alpha=\int(t_1-t_2)\leq 0$. Then
the existence of a metric $H$ on $\cE$ satisfying
\be
i\Lambda F_H=\left(\begin{array}{cc}t_1 \bI_1 & 0\\
	     0&t_2 \bI_2 \end{array}\right),\label{var-extn-eqn}
\ee
implies the $\alpha$-stability of (\ref{extn}).
\end{thm}

It should likewise be possible to adapt the proof of the other
direction of the Hitchin-Kobayashi correspondence.
This will then allow one (by taking $\alpha=0$)  to
establish a more general version of Theorem \ref{existence-tcve}, valid
when $t$\ and $t'$\ are smooth functions satisfying just (\ref{t-t'}).
We will discuss this in a future publication.

To describe the moduli space, let $\cC_E$ and $\cC_L$ the sets of
holomorphic structures on $E$ and $L$ respectively. Consider the set
\be
\cH(E,L)=
\{(\dbar_E,\dbar_L,\phi)\in \cC_E\times\cC_L\times\Omega^0(\Hom(L,E))
\;\;|\;\;
\phi\in H^0(X,\cE\otimes\cL^\ast)\} \label{ht}
\ee
of {\em holomorphic triples} on $(E,L)$, where $\cE$ and $\cL$ denote the
holomorphic vector bundles defined by $\dbar_E$ and $\dbar_L$
respectively.
Let $\cH_\tau(E,L)\subset\cH(E,L)$ be the set of \ts\ holomorphic triples.
This set is invariant under the action of the complex gauge groups
of $E$ and $L$, $\cG^c_E$ and $\cG^c_L$, say. The moduli space of
\ts\  triples is defined as
$$
\cB_\tau(E,L)=\cH_\tau(E,L)/\cG^c_E\times\cG^c_L.
$$

The set $\cB_\tau(E,L)$, which has naturally the  structure of a
variety (cf. \cite{G4}), is closely related to the moduli space
of stable pairs---this is not surprising in view of the definition
\ref{st}.
More precisely, the  map
$(\cE,\cL,\phi)\mapsto (\cE\otimes\cL^\ast,\phi)$ exhibits
$\cB_\tau(E,L)$ as a $\pic^0$-principal bundle over the moduli
space of ($\tau-\deg L$)-stable pairs on $E\otimes L^\ast$.

The study of the general equations (\ref{gcve}), i.e. the case in which
$F$ is of arbitrary rank, requires the introduction of a new notion of
stability. This was carried out in  \cite{BG1}.  All the results
explained above should extend appropriately to the higher rank case when one
replaces $\tau$ and $\tau'$ in (\ref{gcve}) by functions $t$ and $t'$.

%%%%%%%%%%%%%%%%%%%%%%%%%%%%%

%%%%%%%%%%%%%%%%%%%%%%

%\vspace{12pt}

\noindent Department of Mathematics, University of Illinois, Urbana, IL 61801,
USA. bradlow@uiuc.edu

%\vspace{12pt}

\noindent Departamento de Matem\'aticas, Universidad Aut\'onoma de Madrid,
 28049--Madrid, Spain. ogprada@ccuam3.sdi.uam.es

%%%%%%%%%%%%%%

\begin{thebibliography}{ABCD}
%%%%%%%%%%%%%%%%%%%%%%%%%%%%%%
\bibitem[A]{A} M.F. Atiyah, Riemann surfaces and spin structures,
	{\em Ann. Sci. \'Ecole Norm. Sup.} {\bf 4} (1971)
	47--62.


\bibitem[AHS]{AHS}
	M.F. Atiyah, N.J. Hitchin and I.M. Singer,
	Self-duality in four dimensional Riemannian geometry,
	{\em Proc. Royal Soc. Lond.}, Series A, {\bf 362}
	 (1978) 425--61.

\bibitem[Be]{Be}
	A. Bertram, Stable pairs and stable parabolic pairs,
{\em J. Alg. Geom.} {\bf 3} (1994) 703--724.

\bibitem[B1]{B1}
	S.B. Bradlow,
	Vortices in holomorphic line bundles over closed K\"{a}hler
	manifolds,
	{\em Commun. Math. Phys. } {\bf 135} (1990) 1--17.

\bibitem[B2]{B2}
	S.B. Bradlow,
	 Special metrics and stability for holomorphic bundles with
	global sections, {\em \JDG\ } {\bf 33} (1991) 169--214.

\bibitem[BD1]{BD1}
	S.B. Bradlow and G. Daskalopoulos, Moduli of stable pairs for
	holomorphic bundles
       over Riemann surfaces, {\em Int. J. Math.}
       {\bf 2} (1991) 477--513.

\bibitem[BD2]{BD2}
	S.B. Bradlow and G. Daskalopoulos, Moduli of stable pairs for
	holomorphic bundles
	over Riemann surfaces II, {\em Int. J. Math.} {\bf 4} (1993)
       903-925.

\bibitem[BDGW]{BDGW}
	S.B. Bradlow, G. Daskalopoulos,  O. Garc\'{\i}a--Prada
	and R. Wentworth, Stable augmented bundles over Riemann surfaces.
	{\em Vector Bundles in Algebraic Geometry,
	Durham 1993}, Cambridge University Press, 1995.

\bibitem[BG1]{BG1}
	S.B. Bradlow and  O. Garc\'{\i}a--Prada,
	Stable triples, equivariant bundles and dimensional reduction,
	{\em Math. Ann.}, in press.

\bibitem[BG2]{BG2}
	S.B. Bradlow and O. Garc\'{\i}a--Prada,
	Higher cohomology triples and holomorphic extensions,
	{\em Comm. in Analysis and Geom.},  {\bf 3} (1996) 421-464.
\bibitem[Br]{Br}
	R. Brussee,
	$C^\infty$ properties of K\"ahler surfaces,
	 preprint 1995.
\bibitem[D1]{D1}
	S.K. Donaldson,
	Anti-self-dual Yang--Mills connections on a complex algebraic
	      surface and stable vector bundles,
	{\em Proc. Lond. Math. Soc.} {\bf 3} (1985) 1--26.

\bibitem[D2]{D2}
	S.K. Donaldson,
	Infinite determinants, stable bundles and curvature,
	{\em Duke Math. J.} {\bf 54} (1987) 231--247.

\bibitem[D3]{D3}
	S.K. Donaldson,
	The Seiberg--Witten equations and 4-manifold topology,
	{\em Bull. Amer. Math. Soc.} {\bf 33} (1996) 45-70.

\bibitem[FL]{FL}
        P. Feehan and T. Leness,
        Non-abelian monopoles and the relation between Donaldson
        and Seiberg--Witten invariants of smooth four-manifolds,
        in preparation.


\bibitem[FM]{FM}
	R. Friedman and J. Morgan,
	Algebraic surfaces and Seiberg--Witten invariants,
	preprint 1995.

\bibitem[G1]{G1}
	O. Garc\'{\i}a--Prada,
	{\em The Geometry of the Vortex Equation}, D. Phil. Thesis,
	Oxford 1991.

\bibitem[G2]{G2}
	 O. Garc\'{\i}a--Prada,
	 Invariant connections and vortices, {\em Commun. Math. Phys.},
	 {\bf 156} (1993) 527--546.

\bibitem[G3]{G3}
	 O. Garc\'{\i}a--Prada, A direct existence proof for the vortex
	 equations over a compact Riemann surface,
	 {\em Bull. Lond. Math. Soc.} {\bf 26} (1994) 88--96.


\bibitem[G4]{G4}
	 O. Garc\'{\i}a--Prada,
	 Dimensional reduction of stable bundles, vortices and stable
	 pairs, {\em Int. J. Math.}, {\bf 5} (1994) 1--52.

\bibitem[G5]{G5}
	 O. Garc\'{\i}a--Prada,
	 Monopoles and vortices on four-manifolds,
	 {\em Proceedings of the Les Houches summer school
	 on Quantum Symmetries 1995}, Elsevier (Eds. A. Connes and
         K. Gaw\c{e}dzki), to appear.
\bibitem[H]{H}
	N.J. Hitchin,
	 Harmonic spinors,
	{\em Adv. in Math.} {\bf 14} (1974) 1-55.

\bibitem[HL1]{HL1}
	D. Huybrechts and M. Lehn,
	Stable pairs on curves and surfaces, {\em J. Alg. Geom.}
	{\bf 4} (1995) 67--104.

\bibitem[HL2]{HL2}
	D. Huybrechts and M. Lehn,
	Framed modules and their moduli, {\em Int. J. Math.}
	{\bf 6} (1995) 297--324.

\bibitem[JT]{JT}
	A. Jaffe and C. Taubes,
	{\em Vortices and Monopoles}, Progress in Physics {\bf 2},
	Boston, Birkh\"auser, 1980.

\bibitem[KM]{KM}
        P.B. Kronheimer and T.S. Mrowka,
        The genus of embedded surfaces in the projective plane,
        {\em Math. Res. Letts.} {\bf 1} (1994) 797--808.
\bibitem[LM]{LM}
	J.M.F. Labastida and M. Mari\~no,
	Non-abelian monopoles on four manifolds,
	{\em Nuclear Physics} {\bf B 448} (1995) 373.

\bibitem[LaMi]{LaMi}
	H. Blaine Lawson and M.-L. Michelsohn,
	{\em Spin Geometry},
	Princeton University Press, 1989.


\bibitem[OT1]{OT1}
	Ch. Okonek and A. Teleman, The coupled Seiberg--Witten
	equations, vortices, and moduli spaces of stable pairs,
	preprint 1995.

\bibitem[OT2]{OT2}
	Ch. Okonek and A. Teleman, Quaternionic monopoles,
	preprint 1995.

\bibitem[PT]{PT}
	V. Pidstrigach and A. Tyurin,
	Localisation of the Donaldson's invariants along
	Seiberg--Witten classes, preprint 1995.

\bibitem[T1]{T1}
	C.H. Taubes,
	Arbitrary $N$-vortex solutions to the first order
	Ginzburg--Landau equations, {\em Commun. Math. Phys.}
	{\bf 72} (1980) 277--292.

\bibitem[T2]{T2}
	C.H. Taubes,
	On the equivalence of the first and second order equations
	for gauge theories, {\em Commun. Math. Phys.}
	{\bf 75} (1980) 207--227.


\bibitem[Th]{Th}
	M. Thaddeus, Stable pairs, linear systems and the Verlinde formula,
	{\em Invent. Math.} {\bf 117} (1994) 317--353.

\bibitem[UY]{UY}
	K.K. Uhlenbeck and S.T. Yau,
	 On the existence of Hermitian--Yang--Mills connections
	      on stable bundles over compact K\"{a}hler manifolds,
	{\em Comm. Pure and Appl. Math.} {\bf 39--S} (1986) 257--293.

\bibitem[W]{W}
	E. Witten,
	Monopoles and four-manifolds, {\em Math. Res. Letts.}
	{\bf 1} (1994) 769-796.



%%%%%%%%%%%%%%%%%%%%%%
\end{thebibliography}
\end{document}